\RequirePackage[leqno]{amsmath}
\documentclass[leqno]{amsart}
\usepackage{color}
\usepackage{amsfonts}
\usepackage{enumitem}
\setlist[itemize]{topsep=0pt,after=\vspace{1.5\baselineskip}}
 \usepackage[colorlinks=true]{hyperref}
\usepackage{empheq}
\usepackage[T1]{fontenc}
\newtheorem{state}{Assumption}
\def\R{\mathbb R}  
\def
\@cite
#1#2{[{{\bfseries #1}\if@tempswa , #2\fi}]}
\newtheorem{theorem}{Theorem}[section]

\newtheorem{lemma}[theorem]{Lemma}
\newtheorem{proposition}{Proposition}

\newtheorem{remark}{Remark}

\title[On the blow--up time for an attraction--repulsion
chemotaxis system] 
      {Explicit lower bound of blow--up time for an attraction--repulsion
chemotaxis system}
\author[Giuseppe Viglialoro]{}
\RequirePackage[normalem]{ulem} 
\RequirePackage{color}\definecolor{RED}{rgb}{1,0,0}\definecolor{BLUE}{rgb}{0,0,1} 
\providecommand{\DIFaddbegin}{} 
\providecommand{\DIFaddend}{} 
\providecommand{\DIFdelbegin}{} 
\providecommand{\DIFdelend}{} 
\providecommand{\DIFaddbeginFL}{} 
\providecommand{\DIFaddendFL}{} 
\providecommand{\DIFdelbeginFL}{} 
\providecommand{\DIFdelendFL}{} 
\newcommand{\DIFscaledelfig}{0.5}
\RequirePackage{settobox} 
\RequirePackage{letltxmacro} 
\newsavebox{\DIFdelgraphicsbox} 
\newlength{\DIFdelgraphicswidth} 
\newlength{\DIFdelgraphicsheight} 
\LetLtxMacro{\DIFOincludegraphics}{\includegraphics} 
\newcommand{\DIFaddincludegraphics}[2][]{{\color{blue}\fbox{\DIFOincludegraphics[#1]{#2}}}} 
\newcommand{\DIFdelincludegraphics}[2][]{
\sbox{\DIFdelgraphicsbox}{\DIFOincludegraphics[#1]{#2}}
\settoboxwidth{\DIFdelgraphicswidth}{\DIFdelgraphicsbox} 
\settoboxtotalheight{\DIFdelgraphicsheight}{\DIFdelgraphicsbox} 
\scalebox{\DIFscaledelfig}{
\parbox[b]{\DIFdelgraphicswidth}{\usebox{\DIFdelgraphicsbox}\\[-\baselineskip] \rule{\DIFdelgraphicswidth}{0em}}\llap{\resizebox{\DIFdelgraphicswidth}{\DIFdelgraphicsheight}{
\setlength{\unitlength}{\DIFdelgraphicswidth}
\begin{picture}(1,1)
\thicklines\linethickness{2pt} 
{\color[rgb]{1,0,0}\put(0,0){\framebox(1,1){}}}
{\color[rgb]{1,0,0}\put(0,0){\line( 1,1){1}}}
{\color[rgb]{1,0,0}\put(0,1){\line(1,-1){1}}}
\end{picture}
}\hspace*{3pt}}} 
} 
\LetLtxMacro{\DIFOaddbegin}{\DIFaddbegin} 
\LetLtxMacro{\DIFOaddend}{\DIFaddend} 
\LetLtxMacro{\DIFOdelbegin}{\DIFdelbegin} 
\LetLtxMacro{\DIFOdelend}{\DIFdelend} 
\DeclareRobustCommand{\DIFaddbegin}{\DIFOaddbegin \let\includegraphics\DIFaddincludegraphics} 
\DeclareRobustCommand{\DIFaddend}{\DIFOaddend \let\includegraphics\DIFOincludegraphics} 
\DeclareRobustCommand{\DIFdelbegin}{\DIFOdelbegin \let\includegraphics\DIFdelincludegraphics} 
\DeclareRobustCommand{\DIFdelend}{\DIFOaddend \let\includegraphics\DIFOincludegraphics} 
\LetLtxMacro{\DIFOaddbeginFL}{\DIFaddbeginFL} 
\LetLtxMacro{\DIFOaddendFL}{\DIFaddendFL} 
\LetLtxMacro{\DIFOdelbeginFL}{\DIFdelbeginFL} 
\LetLtxMacro{\DIFOdelendFL}{\DIFdelendFL} 
\DeclareRobustCommand{\DIFaddbeginFL}{\DIFOaddbeginFL \let\includegraphics\DIFaddincludegraphics} 
\DeclareRobustCommand{\DIFaddendFL}{\DIFOaddendFL \let\includegraphics\DIFOincludegraphics} 
\DeclareRobustCommand{\DIFdelbeginFL}{\DIFOdelbeginFL \let\includegraphics\DIFdelincludegraphics} 
\DeclareRobustCommand{\DIFdelendFL}{\DIFOaddendFL \let\includegraphics\DIFOincludegraphics} 

\begin{document}
\maketitle
\DIFdelbegin 
\DIFdelend 

\centerline{\scshape Giuseppe  Viglialoro}
\medskip
{\footnotesize
 \centerline{Dipartimento di Matematica e Informatica}
 \centerline{Universit\`{a} di Cagliari}
 \centerline{V. le Merello 92, 09123. Cagliari (Italy)}
  \centerline{E-mail: giuseppe.viglialoro@unica.it}
}

\bigskip
\begin{abstract}
In this paper we study classical solutions to the zero--flux attraction--repulsion chemotaxis--system 
\begin{equation}\label{ProblemAbstract}
\tag{$\Diamond$}
\begin{cases}
u_{ t}=\Delta u -\chi \nabla \cdot (u\nabla v)+\xi \nabla \cdot (u\nabla w) & \textrm{in }\Omega\times (0,t^*), \\
0=\Delta v+\alpha u-\beta v & \textrm{in } \Omega\times (0,t^*),\\
0=\Delta w+\gamma u-\delta w & \textrm{in } \Omega\times (0,t^*),\\
\end{cases}
\end{equation}
where $\Omega$ is a smooth and bounded domain of $\mathbb{R}^2$, $t^*$ is the blow--up time and $\alpha,\beta,\gamma,\delta,\chi,\xi$ are positive real numbers.  
From the literature it is known that under a proper interplay between the above parameters and suitable smallness assumptions on the initial data  $u({\bf x},0)=u_0\in C^0(\bar{\Omega})$, system \eqref{ProblemAbstract} has a unique classical solution which becomes unbounded as $t\nearrow t^*$. The main result of this investigation is to provide an explicit lower bound for $t^*$ estimated in terms of $\int_\Omega u_0^2 d{\bf x}$ and attained by means of well--established techniques based on ordinary differential inequalities.
\end{abstract}
\section{Introduction and motivations}\label{IntroductionSection} 
This paper is dedicated to the following problem
\begin{equation}\label{problem}
\begin{cases}
u_{ t}=\Delta u -\chi \nabla \cdot (u\nabla v)+\xi \nabla \cdot (u\nabla w) & \textrm{in }\Omega, t>0, \\
0=\Delta v+\alpha u-\beta v & \textrm{in } \Omega,  t>0,\\
0=\Delta w+\gamma u-\delta w & \textrm{in } \Omega,  t>0,\\
 u_{\boldsymbol\nu}= v_{\boldsymbol\nu}= w_{\boldsymbol\nu}=0 & \textrm{in } \partial \Omega, t>0, \\
u({\bf x},0)\geq 0 & {\bf x} \in  \Omega,
\end{cases}
\end{equation}
where for the unknown $(u,v,w)=(u({\bf x},t),v({\bf x},t),w({\bf x},t))$ the vector ${\bf x}=(x_1,x_2)$ belongs, unless differently specified, to a bounded and smooth domain $\Omega$ of $\mathbb{R}^2$ and where $\alpha,\beta,\gamma,\delta,\chi,\xi >0$. The function $u_0({\bf x})=u({\bf x},0)$ is nonnegative, sufficiently regular and  corresponds to the initial value of $u$, while the subscription  ${(\cdot)}_{\boldsymbol\nu}$ indicates the outward normal derivative on $\partial \Omega$. 

Like many variants of the well--known models used by Keller and Segel in the celebrated papers \cite{K-S-1970,Keller-1971-MC,Keller-1971-TBC}  to describe general chemotaxis phenomena, system \eqref{problem} represents the situation where the motion in an insulated domain of a certain cell density $u({\bf x}, t)$ at the position ${\bf x}$ and at the time $t$, initially distributed according to the law of $u_0({\bf x})$,  is influenced by the presence of two chemical signal concentrations, $v({\bf x}, t)$ and $w({\bf x}, t)$, which respectively attracts toward the increasing chemoattractant  and repulses from the increasing chemorepellent the same cells. Moreover,  the parameters $\chi$ and $\xi$ measure the strength of the attraction and repulsion, and  the second and third equations idealize that chemoattractant and chemorepellent, $v$ and $w$, are released by cells and  decay with rates $\beta$ and $\delta$. Applications of such a model are met in aggregation phenomena of microglia observed in Alzheimer's disease (see   \cite{Luca2003Alzheimer,PainterHillenCanadianAppMathQuat}).

Strong numerical methods and real experiments indicate that the aforementioned movement may eventually degenerate into aggregation processes, where an uncontrolled gathering of cells at certain spatial locations is perceived as time evolves.  This is the so called \textit{chemotactic collapse}, appearing when $u$, in a particular  instant of time (the blow-up time $t^*$), becomes unbounded in one or more points of its domain. This coalescence phenomena is well understood for the classical parabolic--elliptic Keller--Segel system, obtained by letting $\xi=0$ in \eqref{problem}, and reading in a general bounded domain $\Omega \subset \R^n$, with $n\geq 1$, as:
\begin{equation}\label{problemClassicalK-S}
\begin{cases}
u_{ t}=\Delta u -\chi \nabla \cdot (u\nabla v) & \textrm{in }\Omega, t>0, \\
0=\Delta v+\alpha u-\beta v & \textrm{in } \Omega,  t>0.\\
\end{cases}
\end{equation}
As far as known results tied to this system are concerned,  no blow--up solution can be detected in one--dimensional settings, while in \cite{JaLu} for radial solutions and in \cite{Nagai} for non-radial solutions, the authors prove that for $n=2$ a certain threshold value given by the  product between the chemosensitivity $\chi$ and the initial mass $\int_\Omega u_0 d{\bf x}$ decides whether the solution can blow up at some finite time  or exists for all time $t>0$. 

Unlike \eqref{problemClassicalK-S}, the presence of the attraction--repulsion mechanism in system \eqref{problem} makes the corresponding analysis more complex. To the best of our knowledge, these are the most important achievements obtained up the date, most concerning the planar setting:
\begin{itemize}
\item In high dimensions, precisely $n\geq 2$, if repulsion prevails over attraction, in the sense that $\xi\gamma-\chi\alpha>0$ then for any sufficiently smooth initial data $u_0$, the model possesses globally bounded classical solutions. Conversely, for $n=2$, $\xi\gamma-\chi\alpha<0$ and $\beta=\delta$ there exist appropriate initial data emanating solutions with blow--up at finite time (see \cite{TaoWangAttractionRepulsion});
\item In the bi--dimensional radial case, for $\xi\gamma-\chi\alpha<0$ and any $\beta,\delta>0$ there exist initial data $u_0$ which produce unbounded solutions at finite time (\cite{EspejoSuzukiAttractionRepulsion}). 
\item For the general bi--dimensional case (i.e. removing radial symmetries), $\xi\gamma-\chi\alpha<0$ and any $\beta,\delta>0$ finite time blow--up solutions have been also detected in \cite{YuGuoZhengAttractionRepulsion}. (See also \cite{LiLiAttractionRepulsion}.)
\item In the recent paper \cite{CriticalMassAttrRepulGuoEtAl} it is established that for $n=2$ and $\chi\alpha-\xi\gamma>0$,  the value $\frac{4\pi}{\chi\alpha-\xi\gamma}$ is the critical mass for the attraction--repulsion chemotaxis system \eqref{problem} through which it is possible to identify global boundedness or possible finite time blow--up of solutions.
\end{itemize} 
Motivated by the above  discussion, aim of the present research is expanding the underpinning theory of the mathematical analysis of problem \eqref{problem}, which, so far we are aware, is not included in the above cases.  In particular,  inspired by the presented state of the art, in this work we estimate a lower bound of $t^*$ for classical and unbounded solutions to \eqref{problem}, so to essentially obtain a safe existence interval $[0,t^*)$ where such solutions exist. 

From the technical point of view, starting from those scenarios where local solutions $(u,v,w)$ to \eqref{problem} are detected, we associate to the $u$-component the energy function $E(t):=\int_\Omega u^2 d{\bf x}$ and derive in $(0,t^*)$ a first order differential inequality (ODI); by assuming unboundedness of $E(t)$ in a left neighborhood of $t^*$, an explicit integration will infer the desired lower bound. In the context of estimates of blow--up time to unbounded solutions for evolutive equations, this strategy is rather classical and widely used; in this sense in our computations we will use some well--known ideas and inequalities but, being the general construction of the proof not so straightforward, we take care to make this article self-contained and, further, we are also necessarily required to invoke as many other new derivations and adaptations.

\section{Some premises and preparatory tools }
From  the considerations given above, we continue this paper by presenting the following proposition, whose second part represents the starting point of our work and that we claim according to our purposes. 
\begin{proposition}\label{PropositionExistenceSolution}
Let $\Omega$ be a bounded and smooth domain of $\R^2$. Then, for any $\alpha,\beta,\gamma,\delta\chi,\xi >0$, with $\chi\alpha-\xi\gamma>0$, and nonnegative and nontrivial initial data $u_0({\bf x})\in C^0(\bar{\Omega})$ we have that: 
\begin{enumerate}[label=(i$_\arabic*$)]
\item\label{ItemOneProposition} if  $\int_\Omega u_0({\bf x}) d{\bf x}<\frac{4\pi}{\chi\alpha-\xi\gamma}$, problem \eqref{problem} admits a unique global solution $(u, v,w)\in (C^0(\bar\Omega\times  [0,\infty)) \cap  C^{2,1}(\bar\Omega \times (0,\infty)))^3$ of nonnegative and bounded functions;
\item\label{ItemTwoProposition}   if for some ${\bf x}_{\partial \Omega}\in\partial \Omega$ and appropriate $\rho>0$,  $\int_\Omega u_0({\bf x}) d{\bf x}>\frac{4\pi}{\chi\alpha-\xi\gamma}$ and $\int_\Omega u_0({\bf x}) |{\bf x}-{\bf x}_{\partial \Omega}|^2d{\bf x}<\rho$,  problem \eqref{problem} admits a unique local solution $(u, v,w)\in (C^0(\bar\Omega\times  [0,t^*)) \cap  C^{2,1}(\bar\Omega \times (0,t^*)))^3$ of nonnegative  functions such that for some finite time $t^*$
\begin{equation}\label{BlowUpClassicalSense}
\limsup_{t\rightarrow t^*}\| u(\cdot,t)\|_{L^\infty(\Omega)}=\infty.
\end{equation}
\end{enumerate}
\begin{proof}
These results are, respectively, shown in \cite[Theorem 1 and Theorem 2]{CriticalMassAttrRepulGuoEtAl}.
\end{proof}
\end{proposition}
Having in our hands this existence result, we see that the  interval $I = [0, t^*)$ where classical solutions to system \eqref{problem} are defined can be unbounded ($t^*=\infty$ and the solutions are global, as specified in item \ref{ItemOneProposition} of Proposition \ref{PropositionExistenceSolution})  or bounded ($t^*$ is finite and the solutions are local and blow up as explained in \ref{ItemTwoProposition}). In our contribution, we will study the latter case, precisely by developing an analysis dealing with some estimates for the length of the interval $I = [0, t^*)$. 
\section{Presentation of the main theorem}\label{SectionPresentationResults}
After these considerations, we can present our main result. First, we make the following
	\begin{state}\label{AssumptionDom}
 $\Omega$ is a bounded domain of $\mathbb{R}^2$, star-shaped and convex in two orthogonal directions, whose geometry for some origin ${\bf x}_0$, inside $\Omega$, is defined by 
\begin{equation*}
m_1:= \frac{3}{2 \rho_0},\;\;\; m_2:= 1+ \frac{d}{\rho_0},
\end{equation*}
with $\rho_0 := \min_{\partial \Omega} ( {\bf x}- {\bf x}_0) \cdot {\bf \boldsymbol\nu}$ and $d:= \max_{\overline{\Omega}} | {\bf x}-{\bf x}_0|$.
\end{state}
\begin{theorem}\label{Maintheorem}
Let $\Omega$ be a domain satisfying Assumption \ref{AssumptionDom}. For  $\alpha,\beta,\gamma,\delta,\xi,\chi>0$  such that $\alpha\chi-\xi\gamma>0$ and nonnegative and nontrivial  initial data $u_0({\bf x})\in C^0( \bar \Omega)$ fulfilling item \ref{ItemTwoProposition} of Proposition \ref{PropositionExistenceSolution}, let \[(u, v,w)\in (C^0(\bar\Omega\times  [0,t^*)) \cap  C^{2,1}(\bar\Omega \times (0,t^*)))^3\] be the local solution to system \eqref{problem}, blowing-up at finite time $t^*$ in the sense that
\begin{equation*}
\limsup_{t\rightarrow t^*}\| u(\cdot,t)\|_{L^\infty(\Omega)}=\infty.
\end{equation*}
In such circumstances, if for $E(t):=\int_\Omega u^2 d {\bf x}$ it holds that $\limsup_{t\rightarrow t^*}E(t)=\infty$, then it is possible to find a positive constant $\tilde{c}=\tilde{c}(\delta,\gamma)$ such that
\begin{equation}\label{ExplicitLowerBound} 
t^*\geq\frac{2}{A\sqrt{E(0)}},
\end{equation}
where
\[
A=\Big[\alpha\chi+\frac{\tilde{c} \xi\delta }{3}\Big(\frac{3\gamma}{2\delta}\Big)^{-2}\Big]\frac{\sqrt{2}}{3}m_1+\frac{\tilde{c} \xi\delta }{3}\Big(\frac{3\gamma}{2\delta}\Big)^{-2}.
\]
\end{theorem}
The proof of this theorem involves different general inequalities (see $\S$\ref{PreliminariesSection}), some of these purely associated to properties of functions belonging to specific spaces and of domains where they are defined, others also relying on additional facts, as for instance the type of equation that these functions have to solve. Invoking these relations and making full use of the overall structure of system \eqref{problem}, in $\S$\ref{MainTheoremProofSection} we will derive an energy--type inequality associated to $E(t)$, so that an integration we will enable us to prove, in $\S$\ref{SectionProofsTheorems}, Theorem \ref{Maintheorem}. 
\begin{remark}\label{RemarkOnL2BlowUAndClassical}
As to the connection between the classical blow-up in the $L^\infty(\Omega)$-norm of solutions to \eqref{problem}, i.e. relation \eqref{BlowUpClassicalSense}, and that in the $L^2(\Omega)$-norm (and in general in the $L^p(\Omega)$-norm, $p>1$), i.e. in the sense that $\limsup \int_\Omega u^2 d{\bf x}\nearrow \infty$ as $t\searrow t^*$, we want to observe that once it is assumed that $\Omega$ is a bounded domain, we only know that
\[\| u(\cdot,t)\|_{L^2(\Omega)}\leq |\Omega|^\frac{1}{2}\|u(\cdot,t)\|_{L^\infty(\Omega)},\]
so that if a solution blows up in the $L^2(\Omega)$-norm, it does in the $L^\infty(\Omega)$-norm. Conversely, if a solution becomes unbounded in  the $L^\infty(\Omega)$-norm at some finite time $t^*$, $E(t):=\int_\Omega u^2 d{\bf x}$ may diverge at $t^*$ but might also remain bounded in a neighborhood of it (so to be even continuously prolonged up to the boundary $t^*$). Hence, the lower bound of $t^*$ given in \eqref{ExplicitLowerBound} can be computed under the blow-up assumption in the sense of the $L^2(\Omega)$-norm. In this regard, particular attention should be paid to  how eliminate (or weaken) this extra hypothesis and to how figure out when the blow-up scenario in the $L^\infty(\Omega)$-norm implies that in $L^2(\Omega)$-norm. (Apparently the key to accomplish such purpose 
is an adaptation to our model of a refined extensibility criterion established in \cite[Theorem 2.2]{FREITAGLpLinftyCoincide}.) Nevertheless, since this is not the objective of this paper, we leave herein this question open, and maintain the assumption $\limsup E(t)\nearrow \infty$ as $t\searrow t^*$, exactly following the classical approach used in papers concerning unbounded solutions to general evolutive problems (see, for instance, \cite[Theorem 1 and Theorem 2]{LI-ZhengBlowUpK-S},  \cite[Theorem 2.4 and Theorem 2.7]{MarrasVernierVigliaWithm} and \cite[Theorem 1 and Theorem 2]{PSong} for contributions in the frame of chemotaxis models or \cite[Theorem 2.1]{PAYNEPhilSchaferNonlinAnalysis} and \cite[Theorem 1 and Theorem 4]{PAYNE-Phil-SchaferBounds}  for others in different areas).

\end{remark}
\section{Some functional inequalities: toward the ODI}\label{PreliminariesSection}
We will invoke these two coming lemmas. The first is valid for general functions, with sufficient regularity and defined in suitable domains. On the other hand, 
since in our computations we will chiefly be concerned with the $u$--component, it is desirable to estimate various terms involving the $v$-- and the $w$--components of the solution $(u,v,w)$ to problem \eqref{problem}: this is possible by virtue of the second lemma. 
\begin{lemma}\label{lemma2} 
Let $\Omega$ be a domain satisfying Assumption \ref{AssumptionDom}. Then, for any nonnegative function $V\in C^1(\bar{\Omega})$ we have
\begin{equation}\label{SobolevTypeInequBoundary}
\int_{\partial \Omega} V^2 d s \leq \frac{4m_1}{3}\int_\Omega V^2 d {\bf x}+2(m_2-1)\int_\Omega V|\nabla V|d {\bf x}.
 \end{equation}
Moreover, for any $c_1>0$ it also holds that 
\begin{equation}\label{Inequ_v^3Dim2}
\begin{split}
\int_\Omega V^{3}  d {\bf x} \leq & \frac{\sqrt{2} m_1}{3}\Big(\int_\Omega V^2 d {\bf x} \Big)^\frac{3}{2}+\frac{m_2^2c_1}{16}\Big(\int_\Omega V^2d {\bf x} \Big)^2+\frac{2}{c_1} \int_\Omega |\nabla V|^2 d {\bf x}. 
\end{split}
 \end{equation}
\begin{proof}
As to relation \eqref{SobolevTypeInequBoundary}, we refer to \cite[(A.1)]{PPVPII}.

On the other hand (as in \cite[Lemma 3.2]{VigliaGradTermDiffIntEqua}),  let us consider the inequality following (2.10) in \cite{PS_Neuman_Cond} and, fixing the value of the parameter $n$ therein used equal to 2, we rearrange it as follows: 
 \begin{equation*}
\begin{split}
\Big(\int_\Omega V^{4} d {\bf x}\Big)^\frac{1}{2}\leq & \Big(\frac{1}{2}\int_{\partial \Omega} V^{2}|\nu_1|ds+\int_\Omega V |V_{x_{1}}| d {\bf x}\Big)^\frac{1}{2} 
\times \Big(\frac{1}{2}\int_{\partial \Omega} V^{2}|\nu_2|ds+\int_\Omega V |V_{x_{2}}| d {\bf x}\Big)^\frac{1}{2}.
\end{split}
\end{equation*}
With the identification ${\bf x}=(x_1,x_2)$, applications of the Young inequality give 
\begin{equation}\label{Inequality_u^4_R^2_2}
\begin{split}
&\Big(\int_\Omega V^{4} d {\bf x}\Big)^\frac{1}{2}\leq \frac{1}{4}\Big(\int_{\partial \Omega} V^{2}|\nu_1|ds+\int_{\partial \Omega} V^{2}|\nu_2|ds \Big) 
\\&\quad \quad \quad \hskip 1.3cm + \frac{1}{2}\Big(\int_\Omega V |V_{x_{1}}| d {\bf x}+\int_\Omega V |V_{x_{2}}| d {\bf x}\Big)
\\
&\leq \frac{1}{4} \Big(\int_{\partial \Omega} V^{2} ds\int_{\partial \Omega} V^{2}|\nu_1|^2ds\Big)^\frac{1}{2}+\frac{1}{4}\Big(\int_{\partial \Omega} V^{2} ds\int_{\partial \Omega} V^{2}|\nu_2|^2ds\Big)^\frac{1}{2}\\
&\phantom{=\,\,} +\frac{1}{2} \Big(\int_{\Omega} V^{2} d {\bf x}\int_{\Omega} (V)_{x_{1}}^2d {\bf x}\Big)^\frac{1}{2} +\frac{1}{2}\Big(\int_{\Omega} V^{2} d {\bf x}\int_{\Omega} (V)_{x_{2}}^2d {\bf x}\Big)^\frac{1}{2}\\
&\leq \frac{\sqrt{2}}{4}\int_{\partial \Omega} V^{2} ds+\frac{\sqrt{2}}{2} \Big(\int_\Omega V^{2} d {\bf x}\Big)^\frac{1}{2} \Big(\int_\Omega |\nabla V|^2 d {\bf x}\Big)^\frac{1}{2},
\end{split}
\end{equation}
where in the last step we have also taken into account that $a^\frac{1}{2}+b^\frac{1}{2}\leq \sqrt{2}(a+b)^\frac{1}{2}$ (with $a\geq 0$ and $b \geq 0$).

Now inserting \eqref{Inequality_u^4_R^2_2} in this relation
\[
\int_\Omega V^3 d {\bf x} \leq \Big(\int_\Omega V^4 d {\bf x} \int_\Omega V^{2}d {\bf x} \Big)^\frac{1}{2},
\]
naturally coming from the Hölder inequality, and making use of \eqref{SobolevTypeInequBoundary} leads to 
\begin{equation}\label{SupportEqu}
\begin{split}
\int_\Omega V^3 d {\bf x} \leq&\frac{\sqrt{2}m_1}{3}\Big(\int_\Omega V^2 d {\bf x} \Big)^\frac{3}{2}+\frac{\sqrt{2}(m_2-1)}{2}\Big(\int_\Omega V^2 d {\bf x} \Big)^\frac{1}{2}\int_\Omega V|\nabla V| d {\bf x}\\&+\frac{\sqrt{2} }{2}\int_\Omega V^2 d {\bf x} \Big(\int_\Omega |\nabla V|^2 d {\bf x}\Big)^\frac{1}{2}.
\end{split}
\end{equation}
Finally, since the Hölder inequality infers 
\begin{equation*}
\begin{split}
\Big(\int_\Omega V^2 d {\bf x} \Big)^\frac{1}{2}\int_\Omega V|\nabla V| d {\bf x}&\leq\Big(\int_\Omega V^2 d {\bf x} \Big)^\frac{1}{2}\Big[\Big(\int_\Omega V^2 d {\bf x} \Big)^\frac{1}{2}\Big(\int_\Omega  |\nabla V|^2 d {\bf x} \Big)^\frac{1}{2}\Big] \\&
=\int_\Omega V^2 d {\bf x} \Big(\int_\Omega  |\nabla V|^2 d {\bf x} \Big)^\frac{1}{2},
\end{split}
\end{equation*}
by employing this estimate into \eqref{SupportEqu} and using again Young's inequality with the support of an arbitrarily positive constant $c_1$, we obtain the claimed relation \eqref{Inequ_v^3Dim2}. 
\end{proof}
\end{lemma}
\begin{lemma}\label{EllipticEhrlingSystemLemma}
Let $\delta,\gamma>0$ and $\Omega$ a bounded and smooth domain of $\R^n$, $n\geq 1$. Then there exists $\tilde{c}=\tilde{c}(\delta,\gamma) >0$ such that whenever $f\in C^2(\bar\Omega)$ is nonnegative, the solution $\varphi \in C^{3,\kappa}(\bar\Omega)$, for all $0<\kappa<1$,  of
\begin{equation}\label{EllipticEhrlingSystem}
\begin{cases}
-\Delta \varphi+\delta \varphi=\gamma f & \textrm{in}\quad  \Omega,\\
\varphi_{\boldsymbol\nu}=0 & \textrm{on}\quad  \partial \Omega,
\end{cases}
\end{equation}
satisfies 
\begin{equation*}
\int_\Omega \varphi^3 d{\bf x}\leq \frac{2\gamma^3}{3\delta^{2}}\int_\Omega f^3d{\bf x} +\tilde{c}\Big(\int_\Omega f^2d{\bf x}\Big)^\frac{3}{2}. 
\end{equation*}
\begin{proof}
We reason as in \cite[Lemma 2.2]{WinklerHowFar}. By testing the first equation of \eqref{EllipticEhrlingSystem} by $\varphi^2$, we obtain through the integration by parts formula and the Young inequality
\begin{equation*}
2\int_\Omega \varphi\lvert \nabla \varphi \rvert^2 d {\bf x}+\delta \int_\Omega \varphi^3 d{\bf x}=\gamma \int_\Omega \varphi^2 f d{\bf x}\leq \frac{\delta}{2}\int_\Omega \varphi^3  d{\bf x}+\frac{\gamma^3}{3}\Big(\frac{3\delta}{4}\Big)^{-2}\int_\Omega f^3  d{\bf x},
\end{equation*}
or also
\begin{equation}\label{InequalityLikeEhrling_1}
\int_\Omega \lvert \nabla \varphi^\frac{3}{2} \rvert^2 d {\bf x}\leq \frac{2\gamma^3}{3\delta^{2}}\int_\Omega f^3  d{\bf x}.
\end{equation}
By virtue of the inclusions 
\[W^{1,2}(\Omega)\hookrightarrow \hookrightarrow L^2(\Omega)\hookrightarrow L^\frac{2}{3}(\Omega),\]
Ehrling's Lemma (see \cite[Lemma 1.1]{ShowalterBOOK}) provides for any $\eta>0$ a constant $c(\eta)>0$ such that 
\begin{equation*}
\|V\|_{L^2(\Omega)}^2\leq \eta\|V\|_{W^{1,2}(\Omega)}^2+c(\eta)\|V\|_{L^\frac{2}{3}(\Omega)}^2\quad \textrm{for all } V \in W^{1,2}(\Omega),
\end{equation*}
so that an application of the Hölder inequality explicitly gives for all $V \in W^{1,2}(\Omega)$
\begin{equation}\label{InequalityEhrlinSupport}
\int_\Omega V^2 d{\bf x}\leq \eta\int_\Omega V^2 d{\bf x}+\eta\int_\Omega |\nabla V|^2 d{\bf x}+c(\eta)|\Omega|^\frac{3}{2}\Big(\int_\Omega V^\frac{4}{3}d{\bf x} \Big)^\frac{3}{2}.
\end{equation}
On the other hand, additional standard testing procedures applied again to problem \eqref{EllipticEhrlingSystem} give
\[
\int_\Omega |\nabla \varphi|^2 d{\bf x} +\gamma \int_\Omega  \varphi^2 d{\bf x}= \delta \int_\Omega f\varphi d{\bf x},
\]
and with the support of the Young inequality also
\begin{equation*}
\int_\Omega \varphi^2 d{\bf x} \leq \frac{\gamma^2}{\delta^2}\int_\Omega f^2 d{\bf x}.
\end{equation*}
Finally, using this last relation and \eqref{InequalityEhrlinSupport} with $\eta=\frac{1}{2}$, $\tilde{c}=\tilde{c}(\delta,\gamma)=2 c(\frac{1}{2})|\Omega|^\frac{3}{2}\frac{\gamma^3}{\delta^3}$ and $V=\varphi^\frac{3}{2}$, we can conclude once \eqref{InequalityLikeEhrling_1} is considered.
\end{proof}
\end{lemma}
\subsection{The energy--type ordinary differential inequality}\label{MainTheoremProofSection}
In preparation to the final proof, let us now use all the above derivations to obtain an ODI for the energy function $E(t)=\int_\Omega u^2 d {\bf x}$. This ODI is satisfied by $E(t)$ on the whole $I=[0,t^*)$, both if such energy function is associated to a global solution to system \eqref{problem} than a local; despite this, we will make use of this ODI to derive an explicit estimate for the blow-up time of unbounded solutions. 
\begin{lemma}\label{MainInequalityU^2Lemma}
Let $\Omega$ be a domain satisfying Assumption \ref{AssumptionDom}. Additionally, under the remaining hypothesis of Proposition \ref{PropositionExistenceSolution}, let $(u,v,w)$ be the classical solution of problem \eqref{problem}, with $t^*$ finite or infinite. Then, for $E(t):=\int_\Omega u^2 d{\bf x}$,  $\tilde{c}=\tilde{c}(\delta,\gamma)$ as in Lemma \ref{EllipticEhrlingSystemLemma} and \[
A=\Big[\alpha\chi+\frac{8\gamma \xi\delta }{81}\Big]\frac{\sqrt{2}}{3}m_1+\frac{4 \tilde{c} \xi\delta^3}{27\gamma^2} \textrm{ and }
B=\Big[\alpha\chi+\frac{8\gamma \xi\delta }{81}\Big]^2\frac{m_2^2}{16},
\]
 the following is complied:
\begin{equation}\label{MainInequalityU^2}
\frac{d E(t)}{dt}\leq A E^\frac{3}{2}(t)+BE^2(t)\quad \textrm{ for all } t\in(0,t^*).
\end{equation}
 \begin{proof}
Let us differentiate the functional $E(t)$: we have, using problem \eqref{problem} and the divergence theorem 
\begin{equation}\label{DerivativeE-FirstStep}
\begin{split}
E'(t)&=2\int_\Omega u u_t d {\bf x}=2\int_\Omega u [\Delta u -\chi \nabla \cdot (u\nabla v)+\xi \nabla \cdot (u\nabla w) ]d {\bf x}\\&
=-2\int_\Omega |\nabla u|^2 d{\bf x}+(\alpha\chi-\xi\gamma)\int_\Omega u^3 d {\bf x}+\xi\delta \int_\Omega u^2w d {\bf x}-\chi\beta\int_\Omega u^2 v d{\bf x}\\&
\leq  -2\int_\Omega |\nabla u|^2 d{\bf x}+(\alpha\chi-\xi\gamma)\int_\Omega u^3 d {\bf x}+\xi\delta \int_\Omega u^2w d {\bf x}\quad \textrm{ on } (0,t^*),
\end{split}
\end{equation}
where we neglected the nonpositive term $-\chi\beta\int_\Omega u^2 v d{\bf x}.$ Now, by means of the Young inequality, we have that for any $\varepsilon>0$
\begin{equation*}
\int_\Omega u^2 w d {\bf x}\leq \varepsilon \int_\Omega u^3 d{\bf x}+\frac{1}{3}\Big(\frac{3\varepsilon}{2}\Big)^{-2}\int_\Omega w^3 d{\bf x} \quad \textrm{ on } (0,t^*),
\end{equation*}
so that \eqref{DerivativeE-FirstStep} actually reads 
\begin{equation*}
\begin{split}
E'(t)&\leq  -2\int_\Omega |\nabla u|^2 d{\bf x}+(\alpha\chi-\xi\gamma+\xi\delta\varepsilon)\int_\Omega u^3 d {\bf x}\\ &
\quad + \frac{\xi\delta}{3}\Big(\frac{3\varepsilon}{2}\Big)^{-2}\int_\Omega w^3 d{\bf x}\quad \textrm{ for all } (0,t^*).
\end{split}
\end{equation*}
On the other hand, since the $w$-component solves the third equation of system \eqref{problem}, it is the solution of problem \eqref{EllipticEhrlingSystem} with $f=u$. Hence we estimate the term $\int_\Omega w^3 d{\bf x}$ appearing above  by applying Lemma \ref{EllipticEhrlingSystemLemma}, so to infer 
\begin{equation}\label{DerivativeE-ThirdStep}
\begin{split}
E'(t)&\leq  -2\int_\Omega |\nabla u|^2 d{\bf x}+(\alpha\chi-\xi\gamma+\xi\delta\varepsilon)\int_\Omega u^3 d {\bf x}\\ &
\quad + \frac{2 \xi\gamma^3}{9\delta}\Big(\frac{3\varepsilon}{2}\Big)^{-2}\int_\Omega u^3 d{\bf x}+\frac{\tilde{c} \xi\delta }{3}\Big(\frac{3\varepsilon}{2}\Big)^{-2} \Big(\int_\Omega  u^2 d{\bf x}\Big)^\frac{3}{2}\\ &
\leq  -2\int_\Omega |\nabla u|^2 d{\bf x}+ \tilde{c}_1(\varepsilon)\int_\Omega u^3 d {\bf x}+\tilde{c}_2(\varepsilon) \Big(\int_\Omega  u^2 d{\bf x}\Big)^\frac{3}{2} \quad \textrm{ for all } (0,t^*),
\end{split}
\end{equation}
where
\begin{equation*}
\tilde{c}_1(\varepsilon)=\alpha\chi-\xi\gamma+\xi\delta\varepsilon+\frac{2 \xi\gamma^3}{9\delta}\Big(\frac{3\varepsilon}{2}\Big)^{-2}, \quad \tilde{c}_2(\varepsilon)=\frac{\tilde{c} \xi\delta }{3}\Big(\frac{3\varepsilon}{2}\Big)^{-2}.
\end{equation*}
Finally, if we set $c_1=\tilde{c}_1(\frac{\gamma}{\delta})>0$, we can absorb the addendum involving the gradient in \eqref{DerivativeE-ThirdStep} by invoking  \eqref{Inequ_v^3Dim2} with $V=u$, so arriving at claim \eqref{MainInequalityU^2}.
\end{proof}
\end{lemma}
\section{Proofs of the main result}\label{SectionProofsTheorems}
We are now in the right position to justify our assertion.
\subsubsection*{Proof of Theorem \ref{Maintheorem}} Let $t^*$ be the finite blow-up time of the local solution $(u,v,w)$ to system \eqref{problem}. From Lemma \ref{MainInequalityU^2Lemma}, $u$ satisfies \eqref{MainInequalityU^2} for  any $0<t<t^*$ and, additionally, the assumption $\limsup_{t\rightarrow t^*}E(t)=\infty$ ensures the existence of a time $t_1\in [0,t^{*})$ such that $E(t_1)=E(0)$ and $E(t) \geq E(0)$, for all $t \in [t_1, t^*)$. Subsequently, for all $t\in (t_1,t^*)$ an integration of \eqref{MainInequalityU^2} provides
\begin{equation*}
\int_{t_1}^{t} \frac{d E(\tau)}{A E^\frac{3}{2}(\tau)+BE^2(\tau) }\leq \int_{t_1}^{t} d\tau,
\end{equation*}
or explicitly solving the integrals
\begin{equation*}
\begin{split}
 \frac{2}{A}\Big(\frac{1}{\sqrt{E(t_1)}}-\frac{1}{\sqrt{E(t)}}\Big)+\frac{B}{A^2}\log \Big(\frac{ E(t_1)(A+B\sqrt{E(t)})^2}{E(t)(A+B \sqrt{E(t_1)})^2}\Big)
 \leq t-t_1\leq t.
 \end{split}
\end{equation*}
Finally, taking $t \rightarrow t^*$ and using $\limsup_{t\rightarrow t^*}E(t)=\infty$ and $E(t_1)=E(0)$ we conclude.
 \qed

\subsubsection*{Acknowledgements}
The author is member of the Gruppo Nazionale per l'Analisi Matematica, la Probabilit\`a e le loro Applicazioni (GNAMPA) of the Istituto Na\-zio\-na\-le di Alta Matematica (INdAM) and is partially supported by the research project \textit{Integro-differential Equations and Non-Local Problems}, funded by Fondazione di Sardegna (2017).  
\bibliography{C:/Users/utente/Documents/LavoroUniCa/MIS_ARTICULOS/MIS_Articulos/Bibliography/Bibliography}{}

\end{document}